\documentclass [11pt]{article}
\usepackage{amsmath}
\usepackage{amsfonts}
\usepackage{amssymb}

\newtheorem{theo}{Theorem}
\newtheorem{pro}{Proposition}
\newtheorem{co}{Corollary}
\newtheorem{lema}{Lemma}
\newtheorem{hef}{Definition}
\newtheorem{rem}{Remark}

\begin{document}

\title{Harmonic Analysis of Radon  Filtrations for  $ S_n $ and  $GL_n(q) $ }

\author{M. Francisca Y\'a\~nez\thanks{partially supported by   Fondecyt Grants  1040444, 1070246 and 7070263 and PICS CNRS 1514}}
 \maketitle
\date{}

\begin{abstract}

We present a unified approach to the study of Radon transforms related to symmetric groups and to general linear groups $ \;GL_n(q),\; $ regarded as $q-$analogues of the former.  In both cases, we define a sequence of generalized Radon transforms which are
intertwining  operators for natural representations associated to  Gel'fand spaces for our groups.
This sequence enables us to
decompose in a recursive way these natural representations and to
compute explicitly the associated spherical functions.  Our methods and results   are related by $q-$analogy.
\smallskip


\end{abstract}

\section{Introduction}

In 1917   Radon \cite{ra} showed that a function on   Euclidean space
could be recovered from its integrals over affine hyperplanes,
solving in this way the inversion problem for the nowadays called
``Radon Transform''  which associates to a point function $ \; f
\; $ the hyperplane function $ \; Rf :  H   \mapsto  \int_{H} f. \; $ Later,
Gel'fand \cite{ge} and Helgason \cite{he,hel}
      studied the Radon
Transform in the more general setting of homogeneous spaces for a
Lie group $G$.

In  1975 Soto-Andrade  \cite{jsa}   used   finite Radon Transforms
associated to the geometry of finite symplectic spaces to study the
intertwining algebra of natural representations of the finite similitude  
symplectic group $GSp(4,q)$ in function spaces on isotropic
flags and so obtained by  decomposition the principal series
representations of    $  GSp(4,q).$

In 1977 Dunkl \cite{du2} used finite Radon Transforms to obtain a convenient basis of the intertwining algebra of the natural representation of  $GL_n(q)$ associated to the action of this group on the lattice of subspaces in $n-$dimensional space, to decompose this representation  and to calculate the corresponding spherical functions.

In 1983 Grinberg \cite{gr} independently took up this viewpoint  in the real case to analyze the Radon Transform on compact two
point homogeneous spaces,   considering it   as an
 intertwining operator  between representations of the isometry group
  $ \; G = U(n+1). $
  
Later, in 2002,  Marco and Parcet \cite{marpar}  took advantage of  finite Radon transforms to study the interwining algebra of  the natural representation $L^2( \mathcal P (\Omega))$ of   the symmetric  group  $S (\Omega) $ on a finite set  $ \Omega $, associated to the power set   $\mathcal P (\Omega) $ of   $\Omega$ and to decompose this representation.   

On the other, it should be noticed that  Stanton [8] also considered the Gel'fand spaces $ \; (S_{n},
P_{s}) \; $ and $ \; (GL_n(q), V_s) \; $ and proved that the
spherical functions correspond to discrete orthogonal polynomials, without using Radon Transforms.

We describe now the results in this paper, which were already announced in  \cite{mfy2}. For  $G$ being  the symmetric group or its  $q-$analogue, the finite general linear group,  the Radon transforms that we consider are                                                                                                                                                                                                                                                                                                                                                                                                                                                                                                                                                                                                           $G$ -intertwining operators between natural
representations $ \; \;L^{2}(X) \; \; $ associated to the action
of $G$ on  Gel'fand spaces $ \; X \;\; $ (i.e., $G$ -spaces
whose associated natural representation is multiplicity
free).

We introduce here    Gel'fand spaces $ \; (S_{n}, P_{s}) \;  $ for the symmetric group  $S_{n}$ and their
  ``q-analogues" $ \; \;(GL_{n}(q), V_{s}), \; $ for the finite general linear group $ GL_{n}(q)
$,  where $ \; P_{s} \;
$ is the set of all subsets of $ N = \{ 1,2,...n \}$ with just $
\; s \; $ elements, and $ \; V_{s} \; $ is the set of all
s-dimensional
vector subspaces
of $ \; V = {\bf F}_{q}^{n}, $

\noindent where ${\bf F}_{q}$ is the finite field with $q$ elements, where $q= p^n$, $p$ being a prime number.

We define for each positive integer  $\; s \;$ a ``generalized Radon Transform"  $ \; R_{s}, \; $ between the natural unitary representations $ \; L^{2}
(P_{s}) \; $ and $ \; L^{2} (P_{s+1}) \; $ of $ \; S_{n}. \; $ This is done in an  analogous  way to the construction
of the classical Radon transform.

In theorems 1 and 2 we  prove, applying the  harmonic analysis
techniques   developped in [6], that the sequence of  these
operators and the sequence of their adjoints give a ``resolution''
of the natural representations $ \; L^{2} (P_{s}) \; $ and $ \;
L^{2} (P_{n-s}) \; $, respectively. That is, we show that the
subrepresentations $ \; Ker R_{t}^{*} \; $  (resp.  $ \; Ker R_{n-t}
$) of $ \; S_{n} \; $ give  all the irreducible components of $ \;
L^{2} (P_{s}) \; $ (resp.  $  \; L^{2} (P_{n-s})
 $) . Simultaneously we construct the corresponding spherical
functions. Furthermore, we prove in theorem 3 that the
multiplicity-free natural representations $ \; \; L^{2}(P_{s}) \; \;
$and $ \; \; L^{2}(P_{n-s}) \; \; $ are isomorphic via the operator   $ \mathcal C^*_s  $ induced by taking  complement s. 

In  Theorem  4  we prove mutatis mutandis for  the  $q-$analogue  $\; \;(GL_n(q),V_{s}) \; \; $  of    $ \; \; (S_{n},P_{s}), \; \; $  the results corresponding to those  obtained in Theorems 1 and 2 for $ \; \;(S_{n},P_{s})  $. Finally in  Theorem 5 we prove the $q-$analogue of Theorem 3. However the $q-$analogue    $ \mathcal C^*_s (q) $ of the operator  $ \mathcal C^*_s  $ and  the   proof are  more involved in the case of  $GL_n(q).$    Indeed,    $ \mathcal C^*_s (q) $ is given by a summation over all possible supplementaries  to a given subspace.

\section{The Gel'fand space $ \; (S_n, P_{s}).$}

Let $ \; G = S_{n} \; $ and $ \; N = \{ 1, 2,.... n \}.  $ We
denote by$ \; P_{s} \; $ the set of all subsets of $ \; N \; $
with just $ \; s \; $ elements. We consider the natural,
transitive, action of $ \; S_{n} \; $ on $ \; P_{s}. \; $

We endow $ \; (S_{n}, P_{s}) \; $ with the $ \; S_{n}-$invariant
distance $ \; d\; $ defined by $ \; d(X, X') = |X - X'| = |X|-|X
\cap X'| =\frac{1}{2} |X \triangle X'|). $

\bigskip

\begin{pro} The distance $ d $ is an $ \;  S_{n} - $invariant
function on $\;\; P_{s} \times P_{s}\;\;$ which classifies the
$G-$orbits in $\;\; P_{s} \times P_{s}\;\;$. The number of orbits
is $\;\;s+1\;\;.$

\end{pro}

{\flushleft{\bf Proof:}} It is enough to prove that $ \; (Y, Y') \;
$ and $ \; (Z, Z') \; $ belong to the same orbit iff $ \; |Y \cap
Y'| = |Z \cap Z'|. \; \; \; $ Let $\;\sigma \in
S_{n}\;$ such that $\;\sigma(Y)= Z \;$ and $\;\sigma(Y') = Z'. \;$
Since $ \; \sigma \; $ is bijective, we have $ \; \sigma (Y \cap Y')
= Z \cap Z'. \; $ On the other side, if $ \; |Y \cap Y'| = |Z \cap
Z'| \; $ then we can find a biyection $ \; \sigma : N \rightarrow N
\; $ such that

\[
\sigma (Y-(Y \cap Y'))=Z-(Z \cap Z'); \; \; \; \sigma (Y \cap
Y')=Z \cap Z'; \; \; \; \sigma (Y'-(Y \cap Y'))=Z'-(Z \cap Z')
\]
and $ \; \;\sigma (N - (Y \cap Y')) = N - (Z \cap Z').\; \; $

Then $ \; \sigma (Y) = Z \; $ and $ \; \sigma (Y') = Z'. $

\begin{flushright}     
{$\square  $}

\end{flushright}

\begin{hef}

Let 
$ \; \Omega_{i} (X_{0}) \; $ 
be the set of all elements $ \; X \; $ of
$ \; P_{s} \; $ such that $ \; d(X, X_{0}) = i. \; $

 Furthermore, for $ \;0 \leq t \leq s \leq n \; $ we
define the pseudo-distance $ \; \; \ell \;\;$ on $\;\; P_{t}
\times P_{s}\;\; $  by: $ \; \; \ell(Y,X)=|Y-X| \; \;. $

\end{hef}

\begin{rem}
The number of
elements of $ \; \Omega_{i} (X_{0}) \; $ is $ \; \left( \begin{array}{c} s
\\ s- i \end{array} \right) \left( \begin{array}{c} n-s \\ i \end{array}
\right).$
\end{rem}

\bigskip

\begin{co} The natural unitary representation $ \; (L^{2} (P_{s}), \tau )
\; $ associated to the geometric space $ \; (S_{n} , P_{s}) \; $
is a Gel'fand geometric space, i.e., the unitary representation $
\; ( L^2(P_s), \tau ) \; $, (where $\tau $ is defined by $ \;
(\tau_{\sigma}f)(X)= f(\sigma^{-1}(X)) $ for $ X \in P_s, f \in
L^2(P_s), \sigma \in S_n \; $), is multiplicity free. The number
of irreducible components of
  $\; L^{2} (P_{s}) \; $ is $ \; s+1. $
\end{co}

\bigskip    

Our purpose here is to decompose the natural representation  $ \;
(L^{2} (P_{s}), \tau ) \; $ with the help of certain operators 
which we denote by $ \; R_{s} \; $ and $ \; R^{*}_{s}. \; $ Their
definition is completely  analogous  to the classical Radon
Transform in the continuous case.

\begin{hef}
For $ \; 0 \leq s \leq n, \;  
\; Z \in P_s $ and $ \; X \in P_{s+1}, \; let$

\begin{itemize}

\item[a)] $ R_{s} : L^{2} (P_{s}) \rightarrow L^{2} (P_{s+1}), $
be given by

$ R_{s} (f) (X) = \sum_{Z \subset X} f(Z); \; \; f \in L^{2}(P_{s})
$

and

\item[b)] $  R_{s}^{*} : L^{2} (P_{s+1}) \rightarrow L^{2}
(P_{s}), $ be given by

$ R_{s}^{*} (f) (Z) = \sum_{Z \subset X} f(X); \; \; f \in L^{2}
(P_{s+1}). $
\end{itemize}

\end{hef}

\begin{rem}
 A straightforward calculation
shows that $ \; R_{s} \; $ and $ \; R_{s}^{*} \; $ are
intertwining operators and

\[
R_{s}^{*} \circ R_{s} = (n-s) Id + M_{1}
\]
where $ \; M_{k} = M_{\Omega_{k}} \; $ is the averaging operator
associated to the orbit $ \; \Omega_{k} \; $ of $ \; P_{s}^{2}
/S_{n} \; $ defined by

\[
M_{k} (f)(X) = \sum_{d(X, X')=k} f(X'), \; \; \; (X \in P_{s}).
\]
\end{rem}

\begin{pro} For $ \; 1 \leq s \leq n \; $, we have

\begin{itemize}

\item[i)] $ L^{2} (P_{s}) = Im R_{s-1} \perp Ker R_{s-1}^{*}, $

\item[ii)] $ L^{2} (P_{n-s}) = Im R_{n-s}^{*} \perp Ker R_{n-s}. $
\end{itemize}
\end{pro}

{\flushleft{\bf Proof:}} By computing the Fourier coefficients of a
function $ \; f \in L^{2} (P_{s}) \; $ such that $ \; \langle f, R_{s-1}
\delta_{Z}) \rangle = 0, \; $ where $ \; \delta_{Z} (Z') = \delta_{Z, Z'}
\; $ we obtain

\[
(Im R_{s-1})^{\perp} = Ker R_{s-1}^{*}.
\]
Case ii) is completely  analogous .
\begin{flushright}     
{$\square  $}

\end{flushright}
\bigskip

The following Theorem leads us towards the decomposition of the
multiplicity-free natural representation $\; (L^2(P_m), \tau ) \;
$ into its irreducible components $ \; H^{i} \; (0 \leq i \leq m)
\; $ and gives us the spherical functions $ \; \Phi_{i} \; $
associated to each $ \; H^{i}. $

\bigskip

\begin{theo} Let $ \; 0 \leq m  \leq \frac{n}{2}. \; $ The sequence:

\[
L^{2} (P_{0}) \stackrel{R_{0}}{\rightarrow} L^{2}(P_{1}) \rightarrow
\cdots L^{2} (P_{m-1}) \stackrel{R_{m-1}}{\rightarrow} L^{2}(P_{m})
\]
is an inductive system of monomorphisms which provides all
irreducible components  of $ \; L^{2} (P_{m}) \; $ as successive
orthogonal supplements.
\end{theo}

{\flushleft{\bf Proof:}} In order to prove that $ \; R_{s} \; $ is
injective for $ \; 0 \leq s < \frac{n}{2}, \; $ we will apply
induction. Let $ \; s = 1. \; $ Since $ \; Ker R_{0} \; $ is trivial
we have the decomposition $ \; L^{2} (P_{1}) = H^{o} \perp H^{1}, \;
$ where $ \; H^{o} \; $ is the trivial one dimensional
representation and $ \; H^{1} = Ker R_{0}^{*} \; $ is the Steinberg
representation of dimension $\; n-1\;$ of $ \; S_{n}. \; $We
construct next the spherical functions $ \; \Phi_{0}, \Phi_{1} \; $
associated to $ \; H^{o} \; $ and $ \; H^{1} \; $, respectively. As
we know, $ \; \Phi_{i} \; $ will be the spherical function
associated to $ \; H^{i} \; $ if and only if $ \; \Phi_{i} \; $ is
invariant under the action of the isotropy group $ \; K_{1} \; $ of
the origin $ \; X_{0} = \{ 1 \}. \; $ In our case this means that
the values of $ \; \Phi_{i} \; $ at $ \; X \in P_{1} \; $ depend
only on the distance $ \; d(X_{0}, X). \; $ Furthermore $ \;
\Phi_{i} \; $  should take the values 1 at $ \; X_{0} \; $ and $ \;
\Phi_{i} \;\; $ must belong to $ \; \; H^{i}. \; $ These conditions
allow us to get:

\[
\Phi_{0} \equiv 1 \;
\]
\[
\Phi_{1} (X)  = \left\{ \begin{array}{cc}
1 & X = X_{0} \\
- \frac{1}{n-1} & d(X, X_{0}) = 1
\end{array}
\right.
\]

Next we have to prove that $ \; R_{1} \; $ is an injective operator.
This is however equivalent to proving that $ \; R_{1}^{*} \circ R_{1}  \; $
is an automorphism of $ \; L^{2} (P_{1}). \; $ Since $ \; R_{1}^{*} \circ
R_{1}  = (n-1) Id + M_{1}, \; $ it is enough to prove that the
eigenvalues of $ \; M_{1} \; $ are different from $ \; - (n-1). \; $ We
recall that the eigenvalues $ \; \lambda_{i} \; $ of $ \; M_{1} \; $ are
given by $ \; M_{1} (\Phi_{i}) (X_{0}). \; $ In this way, we obtain: $
\; \lambda_{0} = n-1 \; $ and $ \; \lambda_{1} = 1. \; $ Then $ \; R_{1}
\; $ is injective.

Let us suppose that $ \; R_{t} \; $ is injective for $ \; 0 \leq t <
s; \; s < \frac{n}{2}. \; $ We have to prove that $ \; R_{s} \; $
injective.

 Due to the induction hypothesis and using proposition 1
we have:

\[
L^{2} (P_s) = H^{0} \oplus H^{1} \oplus \cdots \oplus H^{s} \hspace{2cm}
\mbox{(orthogonal sum)}
\]
where $ \;  H^{0} \; $ is the one dimensional trivial representation
and

\[
 H^{t} = (R_{s-1} \circ \cdots \circ R_{t}) \; (Ker
R_{t-1}^{*}), \; (1 \leq t \leq s),
\]

\noindent
 with  dimension $\; \left(
\begin{array}{c} n \\ t
\end{array}\right)-\left( \begin{array}{c} n \\ t-1
\end{array}\right).\;$

We have decomposed $ \; L^{2} (P_{s}) \; $ in $ \; (s+1)-$non
isomorphic subrepresentations of $ \; S_{n}. \; $ Since the number
of $ S_{n}-$orbits in $ \; P_{s}^{2} /S_{n} \; $ is just $ \; s+1,
\; $ we have proved that the subspaces $ \; H^{t}, \; 0 \leq t
\leq s, \; $ are all the irreducible components of $ \; S_{n} \; $
appearing in $ \; L^{2} (P_{s}). \; $  We are now able to
construct the spherical functions $ \; \Phi_{t} \; $ associated to
each $ \; H^{t}. $

We observe that $ \; f \in H^{t} \; $ if and only if there exists
a function $ \; h_t \in Ker R_{t-1}^{*} \; $  such that $ \; f(X)
= \sum_{Y \in P_{t}, Y \subset X} h_t(Y), \; X \in P_{s}. \; $
In order to obtain the spherical function $ \; \Phi_{t} \in
H^{t}, \; $ we need to define a function $ \; h_t\; $ on $ \;
P_{t} \; $ such that $ \; \sum_{Y \in P_{t}, Z \subset Y} h_t(Y) =
0 \; $ for each $ \; Z \in P_{t-1}  \; $ and

\[
\Phi_{t} (X) = \sum_{Y \subset X } h_t(Y) =
\sum_{k=k_{0}}^{min(j,t)} \left( \sum_{Y \subset X, \ell (Y, X_{0})
= k} h_t(Y) \right),
\]
where $ \; X_{0} = \{ 1, 2, \cdots , s \}, \; d(X_{0}, X) = j \; $
and $ \;  \ell (Y, X_{0}) = k_{0} \; $ where $ \; k_{0} = 0 \; $
if $ \; j \leq s-t \; $ or $ \; k_{0} = t+j-s \; $ if $ \; j \geq
s-t. \; $  The values of $ \; \Phi_{t} \; $ at $ \; X \in P_{s} \;
$ must depend only on the  distance$\;d(X_{0}, X).\;$ We impose the
following conditions to the function $ \; h_t: $

\begin{itemize}

\item[i)] $ \sum_{Y \in P_{t}, Z \subset Y} h_t(Y) = 0, \; \;(Z
\in P^{t-1}) $
\item[ii)] $ h_t(Y)= \alpha_t^k \;$ if $\; \ell(Y,X_0) = k \;$

\end{itemize}

Now let us suppose that we have found one such function  $\;h_t
\;$ and let us define the following functions $\;\Phi_{t} \;$ on $\;
P_t \;$ by:

\[
\Phi_{t}(X)= \sum_{Y \subset X}h_t(Y)= \sum
_{k=k_0}^{min(j,t)}|A_{t}^{k,j}|\alpha_t^k ,
\]

\noindent
where

\[
A_{t}^{k,j}= \{Y \in P_t : \ell(Y,X_0)=k, Y \subset X, d(X,X_0)=j.
\]

In order to get an element $ \; Y \in A_{t}^{k,j}, \; $ we need to
choose $ \; t-k \; $ elements among the $ \; s-j \; $ elements of $
\; Y \cap X_{0} \; $ and $ \; k \; $ elements among the $ \; j \; $
elements of $ \; Y - X_{0}. \; $ Therefore  $ \; |A_{t}^{k,j}| =
\left( \begin{array}{c} s-j \\ t-k
\end{array}\right) \left( \begin{array}{c} j \\ k \end{array} \right).  $
\\

Since for all $ \; \sigma \; $ belonging to the isotropy group $ \;
K_{t}, \; $ we have $ \;  \; d(\sigma X, \sigma X_{0}) = d(X, X_{0})
\; $ and $ \; |A_{t}^{k,j}| \; $ depends only on the number of
elements of $ \; X \; $ and $ \; Y, \; $ we obtain that $ \; {
\Phi_t} \; $ is invariant under the action of $ \; K_t. \; $ Then,
for $ \; \Phi_t \; $to be a spherical function, we need still only
to fulfill the normalization condition $ \; \Phi_{t} (X_{0}) = 1 \;
$ or equivalently:

\begin{itemize}

\item[iii)] $ A_{t}^{0,0} \alpha_{t}^{0} = 1. $
\end{itemize}

So we obtain the function $ \; h_t \; $ by solving the following
system of linear equations arising from i), ii) and iii)

\[ \left\{
\begin{array}{l}
|C_{t}^{k,k}| \alpha_{t}^{k} +|C_{t}^{k+1,k}| \alpha_{t}^{k+1} = 0,  \hspace{3cm} 0 \leq
k
\leq t-1\\
\\
\alpha_{0} = \left( \begin{array}{c} s \\ t \end{array} \right)^{-1},
\end{array}
\right.
\]
where

\[
C_{t}^{k',k} =  \{ Y \in P_{t}: Z \subset Y, \; \ell (Y, X_{0}) = k' , \;
\ell (Z, X_{0}) = k \}, \; Z \in P_{t-1}.
\]

  We have $ \; Y = Z \cup \{ y \}. \; $ Then if $ \;y \in X_{0} \; $ we get $ \; k' = k, \; $ and if $ \;y \notin X_{0} \; $ then $ \; k' = k+1. \; $ A
combinatorial computation gives us that
 $ \; |C_{t}^{k,k}| = \left(
\begin{array}{c} s-(t-1-k) \\ 1 \end{array} \right) \; $
and $ \; |C_{t}^{k+1,k}| = \left(
\begin{array}{c} n-(s+k) \\ 1 \end{array} \right). $

Solving this system we get

\[
\alpha_{t}^{0} = \left( \begin{array}{c} s \\ t \end{array} \right)^{-1}
\; \; \; \mbox{and} \; \; \;
\alpha_{t}^{k} = \frac{(-1)^{k} (s-t+k)!(n-s-k)!}{\left( \begin{array}{c}
s \\ t \end{array} \right) (s-t)! (n-s)!}, \; \; \; 1 \leq k \leq t.
\]

Then the spherical function $ \; \Phi_{t} \; $ associated to $ \;
H^{t} \; $ is given by

\[
\Phi_{t} (X) = \sum_{k=k_{0}}^{min(j,t)} (-1)^{k} \left(
\begin{array}{c} s-j \\ t-k
\end{array} \right) \left( \begin{array}{c} j \\ k \end{array} \right)\gamma_{t}^{k}
\]

where $ \; j = d(X_{0}, X ) \;$ and $\;\gamma_{t}^{k} = \left( \begin{array}{c} s \\ t \end{array} \right)^{-1}
\frac{(s-t+k)!(n-s-k)!}{(s-t)!(n-s)!}. $

\bigskip

We have to compute now the eigenvalues of $ \; M_{1} \; $ with
respect to our decomposition of $ \; L^{2} (P_{s}) . \; $ As we
know, the eigenvalues $ \lambda _t $ associated with $ \; H^{t} \; $
are given by $ \; \lambda_{t} = M_{1} (\Phi_{t})(X_{0}). \; $ Then
we have:

\[
\lambda_{t} = (n-s)(s-t) - (s-t+1)t, \; \; \; \; 0 \leq t \leq s
\]

We want to show that $ \; \lambda_{t} \neq -(n-s). \; $ Let us
suppose that for some $ \; t \; $ the eigenvalue $ \; \lambda_{t}
= -(n-s). \; $ By replacing $ \; \lambda_{t} = - (n-s) \; $ in the
last equation and dividing by $ \; s-t+1 \neq 0, \; $ we obtain
that $ \; t+s = n. \; $ But since by hypothesis $ \; t+s < n, \; $
we conclude that $ \; \lambda_{t} \neq - (n-s); \; 0 \leq t \leq s
\; $ and $ \; R_{s}^{*} \circ R_{s} \; $ is an automorphism. Then
$ \; R_{s} \; $ is injective  for all $ \; 0 \leq s < m \; $ and

\[
L^{2} (P_{m}) \cong (Im R_{0}) \oplus (Im R_{0})^{\perp } \oplus
\cdots  \oplus (Im R_{m-1})^{\perp}\;\;\;  (orthogonal\; sum)
\]
       
In this way, for each $ \; 0 \leq m \leq \frac{n}{2}, \; $ we have
constructed all the irreducible components for the natural
representation $ \; L^2(P_{m}) \; $ of$ \; S_n \; $  and the
associated spherical functions.

\begin{flushright}     
{$\square  $}

\end{flushright}

\vspace {12pt} In order to prove that for $ \; 0 \leq m \leq
\frac{n}{2} \; $

\[
L^{2} ( N ) \stackrel{R_{n-1}^{*}}{\rightarrow} L^{2} (P_{n-1})
\rightarrow \cdots \stackrel{R_{n-m}^{*}}{\rightarrow} L^{2}
(P_{n-m}),\;
\]
is a resolution for the Gel'fand space $ \; L^{2}
(P_{n-m}) \;, $ we consider the following preliminaries

\vspace{12pt}

\begin{pro} Let us define the $ \; S_{n}-$invariant
application $ \; {\cal C}_s \; $ from $ \; P_{s} \; $ to $ \; P_{n-s}, \;
(0 \leq s \leq \frac{n}{2}) \; $ by $ \; {\cal C}_s (X) =
N-X. \; $  The following combinatorial properties of $ \;
{\cal C}_s \; $ hold:

Let $ \; Y \in P_{t}, \; X_{0}, X \in P_{s} \; $ and $ \; Z \in P_{t-1} $

\begin{enumerate}

\item If $ \; \ell (Y, X_{0}) = k \; $ then $ \; \ell({\cal C}_s (X_{0}),
{\cal C}_s(Y)) = k \; \; $

\item Let $ \; B_{t}^{k,j} = \{ Y' \in P_{n-t}: Y' \supset X' , \; d(X',
X'_{0}) = j, \; \ell (X'_{0},Y') = k \} \; $ where $ \; X', X'_{0}
\in P_{n-s}, \; $  then $ \; |A_{t}^{k,j}|=|B_{t}^{k,j}|. \; \;$

\item Let $ \; D_{t}^{k',k} = \{ Y' \in P_{n-t}: Y' \subset Z', \;
\ell(X'_{0},Z') = k, \; \ell (X'_{0},Y') = k' \} \; $ where $ \; Z'
\in P_{n-(t-1)},\; X'_{0} \in P_{n-s}), \; \;$then  $ \;
\;|C_{t}^{k',k}|= |D_{t}^{k',k}|,\; \;k'=k,k+1. \; \;$

\item The map $ \;{\cal C}_s \; $ is an $ \; S_{n}-$invariant bijection and $ \;{\cal C}_s^{-1}= {\cal C}_{n-s}\; $.
\end{enumerate}
\end{pro}

\vspace {12pt}

{\flushleft{\bf Proof:}} Since $ \; \; (N-X_{0})-(N-Y)=Y-X_{0} \;
\;$ we get the first statement. In connection with the second
statement, we note that, in order to find $ \; \;Y' \in B_{t}^{k,j},
\; \; $we need to choose $ \; \;j-k\; \;$elements of the $\; \; j \;
\; $elements of $ \; \;X'_{0}-X' \; \;$and we need to select the $
\; \; (n-t)-(n-s+j-k)\; \; $remaining elements of the $\;
\;n-(n-s+j)\; \; $ elements of $ \; \;N-(X'_{0} \cup (X'-X'_{0})).\;
\;$Then $ \; \; |B_{t}^{k,j}| = \left( \begin{array}{c} s-j
\\s-j-(t-k)
\end{array}\right) \left( \begin{array}{c} j \\j-k \end{array} \right)=
|A_{t}^{k,j}|. \;\; $

\bigskip

Finally, if $ \; \; Y'\in D_{t}^{k',k} \; \; $ then $ \; \; Y'=Z'-\{z\} \;
\; $ and $
\; \; (X'_{0}-Y')=(X'_{0}-Z') \cup (X'_{0} \cap \{z\}). \; \; $ Now, if
$ \; \; z \in X'_{0} \; \; $ then $ \; \; |X'_{0}-Y'|=k+1 \; \; $ and if $
\;
\; z \notin X'_{0} \; \; $then $ \; \; |X'_{0}-Y'|=|X'_{0}-Z'|=k. \; \; $
Therefore $ \; \; |D_{t}^{k,k}|= n-(t-1)-(n-s-j)=s-t+j+1 \; \; $and$ \; \;
|D_{t}^{k+1,k}|=(n-s-j). \; \; $

\begin{flushright}     
{$\square  $}

\end{flushright}

\vspace {12pt}

Using these properties and following the same procedure employed in the
proof of Theorem 1 we get the following Theorem:

\vspace {12pt}

\begin{theo} Let $ \; 0 \leq s \leq \frac{n}{2} $

\begin{enumerate}

\item $ R_{n-s}^{*} \; $ is an injective intertwining operator.

\item $ L^{2} (P_{n-s}) \cong Im R_{n-1}^{*} \oplus Ker  R_{n-1} \oplus \cdots
\oplus Ker R_{n-s}. $ (orthogonal sum)

\item The spherical function $ \; \varphi_{t} \; $ associated to
each irreducible sub-representation $ \; L^{t} \cong Ker R_{n-t} \;
$ of $ \; (L^{2} (P_{n-s} (N), S_{n}) \; $ is given by $ \;
\varphi_{t}(X') =  \Phi_{t}({\cal C}_{n-s}(X') , \;
  \; X' \in P_{n-s}. \; \; $

\end{enumerate}
\end{theo}
\begin{flushright}     
{$\square  $}

\end{flushright}

\vspace{12pt}

\begin{theo}    Let $ \; 0 \leq s \leq \frac{n}{2} \; $
$ \; L^{2} (P_{s}) \cong L^{2} (P_{n-s}) $ and $ \; Ker R_{n-s} \cong Ker R_{s-1}^{*}. $
\end{theo}

\vspace{12pt}

{\flushleft{\bf Proof:}} We define $ \; {\cal C}_{s}^{*} : L^{2} (P_{n-s})
\rightarrow L^{2} (P_{s}) , \; 0 \leq s \leq \frac{n}{2}\; $ by $ \;
{\cal C}_{s}^{*}(f)=f \circ {\cal C}_{s}, \; (f \in L^{2} (P_{n-s}). \; $ Let $ \; h \in Ker
R_{n-s} \; $ and $ \; Z \in P_{s-1}\; \; $ then

\[
\sum_{Z \subset Y} {\cal C}_{s}^{*} (h) (Y) = \sum_{Z \subset Y} h(N-Y) =
\sum_{Y' \subset N-Z} h(Y') = 0.
\]

Therefore $ \; {\cal C}_{s}^{*} (Ker R_{n-s}) = Ker R_{s-1}^{*}. \;
$ From this and from Theorems 1 and 2 we get this Theorem.

\begin{flushright}     
{$\square  $}

\end{flushright}

\section{The q-analogue Gel'fand space $ (GL_n(q), V_{m}).$}

Let $ \; V = {\bf F}_{q}^{n} \; $ and $ \; V_{s} \; $ the set of all
$ s-$dimensional subspaces of $ \; V \;, 0 \leq s \leq n. \; $ We
fix a basis $ \{e_1, e_2,...., e_n \}\; $ of $\;V \; $ and let $ \;
W_0 \in V_s,\;$ be the subspace of $ \; V \; $ generated by $ \{e_1,
e_2,....,e_s \} \;$ and $\; W'_0 \in V_{n-s} \; $ the subspace
generated by $ \; \{ e_{s+1},...., e_n \}. \; $

\bigskip

We construct a resolution for the Gel'fand space $ \; L^{2}
(V_{m}, GL_{n} (q)), \; $ through the generalized Radon transform
following the same procedure used for $ \; (L^{2} (P_{m}), S_{n}),0 \leq m \leq \frac{n}{2}.
\; $ The combinatorial calculations are obtained by taking the
q-analogue of $ \; \left(
\begin{array}{c} n \\ m \end{array} \right). \; $

\vspace{12pt}

\begin{hef}

\begin{enumerate}

\item We define the q-analogue of $ \;n! \;$  by $ \; n!_{q} = n_{q}(n-1)_{q} \cdots 1_{q}  \;$
where $ \; n_{q} = \frac{q^{n}-1}{q-1} \;$ and the q-analogue of $\;\left( \begin{array}{c}
n \\ m \end{array} \right) \:$ by

\[
\left( \begin{array}{c}
n \\ m \end{array} \right)_{q} = \frac{n!_{q}}{(n-s)!_{q} s!_{q}}
\]

\item Let $ \; \; (0\leq s \leq m).\;\;$ We define  the q-analogue
pseudo distances $ \; \ell_{q} \;  $ on $ \; V_{s} \times V_{m},
\; \;\; \; $ by $ \;\ell_{q} (U, W) = dim U - dim( U\cap W), \;
\;\;(U \in V_{s}, W \in V_{m}). \;\;\; $

\end{enumerate}
\end{hef}

\vspace{12pt}
\begin{rem}
%{\flushleft{\bf Remarks:}}
\begin{enumerate}

\item The gaussian binomial coefficient $ \;\left(
\begin{array}{c} n \\ s \end{array} \right)_{q} \; $ gives the
number of $s-$dimensional subspaces of $\; V. \; $

\item If $ \; s = m \; $ then $ \; \ell_{q} = d_{q} \; $ defines a
distance on $ \; V_{m} \times V_{m} \; $ which classifies the
orbits of the action on $ \; V_{m} \times V_{m} \; $ of $ \;
GL_n(q). \; $

\item  The natural unitary representation $ \; (L^{2} (V_{s}), \tau )
\; $ associated to the geometric space $ \; (GL_n(q) , V_{s}) \; $
is a Gel'fand geometric space, i.e., the unitary representation $
\; ( L^2(V_s), \tau ) \; $, (where $\tau $ is defined by $ \;
(\tau_{g}f)(W)= f(g^{-1}(W)) $ for $ W \in V_s, f \in
L^2(V_s), g \in GL_n(q) \; $), is multiplicity free. The number
of irreducible components of
  $\; L^{2} (V_{s}) \; $ is $ \; s+1. $

\end{enumerate}
\end{rem}
Let us consider the following Lemma.

\begin{lema}  Let $ \; W \in V_{s}. \; $ The number of $ \; l-$dimensional
subspaces $ \; Z \; $ of $ \; V \; $ such that the dimension of $ \; Z\cap W \; $ is $ \; t \; $ is given by:

\[
q^{(s-t)(l-t)} \; \left( \begin{array}{c} n-s \\ l - t
\end{array} \right)_{q} \left( \begin{array}{c} s \\ t \end{array}
\right)_{q}.
\]
\end{lema}

{\flushleft{\bf Proof:}}

Consider a $\;t-$dimensional subspace $\; U \;$ of $ \;  W. \; $
There are  $\;\; \left( \begin{array}{c} s \\ t \end{array}
\right)_{q}\;\;$
 such 
  subspaces. Let us choose the basis
$\;\{u_1,..,u_t,w_1,..,w_{s-t},v_1,..,v_{n-s}\} \;$ of
$\; V, \;$ where $\; U = <u_1,..,u_t> \;$ and $\;W = U \oplus
<w_1,..,w_{s-t}>. \;$ In order to 
 construct 
$ l-$dimensional
subspaces $\; Z \;$ so that $\; U = Z \cap W, \; $  we need to add
$\;l-t\;$ linearly independent vectors, $\;z_1,..,z_{l-t},\;$ to
$\; \{u_1,..,u_t \}\;$ satisfying
 $\; dim(<z_1,..,z_{l-t}> \cap W )= 0. \;$

\bigskip

Let
 $\;z_i = {\sum_{k=1}^{t}\alpha_k^{i}u_k} + {\sum_{k=1}^{s-t}\beta_k^{i}w_k} + {\sum_{k=1}^{n-s}\gamma_k^{i}v_k}.\;$
Since  $\;z_i \notin <w_1,..,w_{s-t}>,  \;$ then $\;\;(
\alpha_1^i,..,\alpha_t^i,\gamma_1^i,..,\gamma_{n-s}^i ) \neq
\vec{0}, 1 \leq i \leq (l-t).\;\;$ Moreover, if we redefine $\; z'_i
= z_i - {\sum_{k=1}^{t}\alpha_k^{i}u_k} \;$ it is verified that
$\;<u_1,..,u_{t},z_1,..,z_{l-t}> =
<u_1,..,u_{t},z'_1,..,z'_{l-t}>. \;$ Thus there are

\bigskip

$$ \frac{q^{s-t}(q^{n-s}-1)q^{s-t}(q^{n-s}-q)\cdots,q^{s-t}(q^{n-s}-q^{l-t-1})}{(q^{l-t}-1)\cdots (q^{l-t}-q^{l-t-1})}
$$
\noindent

\bigskip
\noindent
ways to complete $\;\;U, \;\;$ from which the lemma follows.
\begin{flushright}     
{$\square  $}

\end{flushright}

\vspace{12pt}

\begin{co} Let $ \;W \in V_{s},\;\; W = W_1 \oplus W_2.\;\;$ The number of $k-$dimensional
subspaces $ \; Z \; $ of $\; W \;$ such that $\; \;\ell_{q} (Z, W_1)=k \;\;$ is given by

\[
q^{(dim W_1)k} \; \left( \begin{array}{c} dim W_2 \\ k
\end{array} \right)_{q}.
\]
\end{co}

\vspace{12pt}

\begin{hef}
 Next we consider the following sets,
( which are q-analogue to the sets $ \; \; \Omega_{i}(X_{0}),
 \; A_{t}^{k,j},\;C_{t}^{k,k'}, \;B_{t}^{k,j},\;D_{t}^{k,k'} \;\; $
defined in the above section).

\begin{enumerate}

\item $\Omega_{i}^{q} (W_{0})= \{ W \in V_{s} | d_{q} (W, W_{0}) = i \}$

\item $A(q)_{t}^{k,j} = \{ U \in V_{t} | \ell_{q} (U, W_{0}) = k, \; U
\subset
W, \; d_{q} (W, W_{0}) = j \}$

where  $ \;\; W \in V_{s}. \;\;$

\item $C(q)_{t}^{k',k} = \{ U \in V_{t} | Z < U, \ell_{q} (U, W_{0}) =
k', \; \ell_{q} (Z,  W_{0}) = k \}$

where $ \; Z \in V_{t-1} \; $ and $ \; k' = k, \; k + 1 $

\item $B(q)_{t}^{k,j} = | \{ U' \in V_{n-t} | W' < U', d_{q} (W',W'_{0}) = j,
\ell_{q} (W'_{0}, U') = k \}$

where  $ \;\; W' \in V_{n-s}. \;\;$

\item $D(q)_{t}^{k',k} =  \{ U' \in V_{n-t} | U' < Z', \ell_{q} (W'_{0},
Z') = k, \ell_{q} (W'_{0}, U') = k' \}$

where $ \; Z' \in V^{n-(t+1)} \; $ and $ \; k'=k+1, \; k. $

\end{enumerate}
\end{hef}
\vspace{12pt}

The following proposition gives us the number of elements of the $q-$analogue sets defined above,
getting similar results to those  contained in proposition 3.

\begin{pro}

\begin{enumerate}

\item $ |\Omega_{j}^{q} (W_{0})| = \left( \begin{array}{c} s \\ s-j
\end{array} \right)_{q} \left( \begin{array}{c} n-s \\ j\end{array}
\right)_{q}  q^{j^{2}} $

\item $ | A(q)_{t}^{k,j}| = \left( \begin{array}{c} s-j \\t-k \end{array}
\right)_{q} \left( \begin{array}{c} j \\ k \end{array} \right)_{q}
q^{((s-j)-(t-k))k} $

\item $ | C(q)_{t}^{k,k}| = \left( \begin{array}{c} s-(t-1-k) \\ 1 \end{array}
\right)_{q}  $

\item $ | C(q)_{t}^{k+1,k}| = \left( \begin{array}{c} n-(s+k) \\ 1 \end{array}
\right)_{q}  q^{s-(t-(k+1))} $

\item If $ \; 0 \leq s \leq n/2, \; $ then $ \; |A(q)_{t}^{k,j}| = |
B(q)_{k}^{j}|, \; \;
|D(q)_{t}^{k,k}| = |C(q)_{t}^{k,k}| \; $ and $ \; |D(q)_{t}^{k+1,k}| =
|C(q)_{t}^{k+1,k}|. $

\end{enumerate}
\end{pro}
\bigskip

{\flushleft{\bf Proof:}} Let $\; V = W_0 \oplus W'_0. \;$ We have
that $\;dim(W \cap W_0) = s-j. \;$ Applying the Lemma, for $\; l=s
\;$ and $\;t=s-j, \;$ we get the first item.

\bigskip

If $\; U \in A(q)_{t}^{k,j} \;$ then $\; U \cap W_0 < W \cap W_0. \;$ Let $\;X \;$
be a $k-$dimensional subspace of  $\; W \cap W_0. \;$ In order to get a subspace $\; U \;$
of $\; A(q)_{t}^{k,j}, \;$ we have to complete each subspace  $\;\;X = U \cap W_0 \;\;$
with a supplement $\; X' \;$ so that $\; X'< W \;$ and
$\; X' \cap W_0 = \{ \vec{0} \}.  \;$ Let $\; W \cap W_0 = X \oplus Y, \;$ then
$\; W = X \oplus Y \oplus Z. \;$ So, we have that $\; X' \;$ is a subspace of
 $\;  Y \oplus Z, \;$ such that $\; X'\cap Y = \{ \vec{0} \}. \;$
Applying the above lemma we get the second item.

\bigskip

To prove item 3, we remark that in order to 
 construct
 a subspace  $\; U \in C(q)_{t}^{k,k} \;$
we need to complete  $\; Z \;$ with a one dimensional subspace
$\; < v > \;$ such that $\; v \in W_0 \;$ but
$\; v \notin (Z \cap W_0). \;$ Let $\; W_0 = (Z \cap W_0) \oplus Y. \;$ Since  $\; v \notin Z,\;$ in order  to complete $\; Z. \;$ it is enough to consider all the one dimensional
subspaces $\; < v > \;$ of $\; Y. \;$ From this 
  item 3  
 follows.

\bigskip

To prove item 4 let $\; Z = (Z \cap W_0) \oplus X ,\;$ and
$\; W_0 = (Z \cap W_0 ) \oplus Y .\;$ Since $\; X \cap W_0 = \{ \vec{0} \} \;$
we have that $\; V = Z \oplus Y \oplus Y'.\;$ To get $\; U \;$ we need that
 $\;v \notin (Z \oplus Y).\;$ Moreover, since $\;Z < U \;$
it is 
 enough 
  that $\; v = y + y' \;$ where $\; y \in Y, \;$
$\; y'\in Y' \;$ and $\; y'\neq \vec{0}. \;$ Then item 4 follows.

\bigskip

 analogous ly, we calculate that:
$$ | B(q)_{t}^{k,j}| = \left( \begin{array}{c} s-j \\s-j-(t-k) \end{array}
\right)_{q} \left( \begin{array}{c} j \\ j-k \end{array} \right)_{q}
q^{((s-j)-(t-k))k}, $$

\bigskip

$$ | D(q)_{t}^{k,k}| = \left( \begin{array}{c} s-(t-1-k) \\ s-t-k \end{array} \right)_{q} $$

 and

\bigskip

$$\:\; | D(q)_{t}^{k+1,k}| = \left( \begin{array}{c} n-(s+k) \\ 1 \end{array}
\right)_{q}  q^{s-(t-(k+1))}. $$

getting item 5.

\begin{flushright}     
{$\square  $}

\end{flushright}

\vspace{12pt}

{\flushleft{\bf 2.1.3. q-Radon 
Transforms 
$ \;  R(q)_{s} \; $ and $
\;
R(q)_{s}^{*} .$}} \\

We define $ \; R(q)_{s} \; $ and $ \; R(q)_{s}^{*} \; $ in an
 analogous  way to $ \; R_{s} \; $ and $ \; R_{s}^{*}. \; $ We
replace $ \; P_{s} \; $ by $ \; V_{s} \; $ and $ \; X \subset Y \;
$ by $ \; U < W \; $ (subspace of $ W$). As we can expect, similar
properties to those contained in Propositions 1 and 2 are
satisfied for $ \; R(q)_{s} \; $ and $ \; R(q)_{s}^{*}. \; $
Following the same procedure used to prove Theorem 1 and using
Proposition 4 and Corollary 2 we prove the following Theorem:

\vspace{12pt}

\begin{theo} Let $ \; 0 \leq s \leq n/2 ; $

\begin{itemize}
\item[a)] If $\; H^{t}(q) = ( R(q)_{s-1} \circ
\dots \circ R(q)_t)(Ker R(q)_{t-1}^{*}\;$ for $\; 1 \leq t \leq s  \; $ then

\[
L^{2} (V_{s}) \cong Im R(q)_{0} \oplus Ker  R(q)_{0}^{*} \oplus \cdots
\oplus Ker R(q)_{s-1}^{*}\;\;\;\;(orthogonal\; sum)
\]
and the spherical funtion  $ \; \Phi(q)_t   \;$ asociated with $\; H^{t}(q)\cong Ker R(q)_{t-1}^{*} \;$
is  given by

\[
\Phi(q)_{t}(W) = \sum_{k=k_{0}}^{min(t,j)} (-1)^{k} q^{\frac{k}{2}(k-1-2j)} \left(
\begin{array}{c} s-j \\ t-k \end{array} \right)_{q}  \left(\begin{array}{c} j \\ k \end{array} \right)_{q}\gamma(q)_{t}^{k}
\]

where $ \; W \in V_s ,  d (W_{0}, W) = j, \;\;\;\gamma(q)_{t}^{k} = \left(\begin{array}{c} s \\ t \end{array} \right)_{q}^{-1}
\frac{(s-t+k)!_{q}(n-s-k)!_{q} }{(s-t)!_{q} (n-s)!_{q}}
 $ and $ \; \; k_{0} = 0
\; \; $if$\; \; t \leq (s-j) \; \; $ and $\; \; k_{0}= t+j-s \; \;
$ if $ \; \;  t \geq (s-j). \; \; $

\item[b)] If $\; H^{n-t}(q) = ( R(q)_{n-s}^{*} \circ
\dots \circ R(q)_{n-t-1}^{*})(Ker R(q)_{n-t})\;$ for $\;1 \leq t \leq s  \; $ then
\[
L^{2} (V_{n-s}) \cong Im R(q)_{n-1}^{*} \oplus Ker  R(q)_{n-1} \oplus \cdots
\oplus Ker R(q)_{n-s,}
\]

and the spherical function $ \; \varphi(q)_t \;$ asociated with $\; H^{n-t}(q)\cong Ker R(q)_{n-t}\;$, is given by

\[
\varphi(q)_{t}(W') = \Phi(q)_{t}({W})
\]

\noindent
where:$ \; \; W \in V^{s} \; \;$ and $  \; \; d_{q}(W',W'_{0})=
d_q(W,W_{0})=j.\; \; $

\end{itemize}
\end{theo}

{\flushleft{\bf Proof:}} Following the same procedure used in theorem 1 to get that if $\; R_t\;$ is injective for $\;1 \leq t < s  \; $then $\;  L^2(P_s)= H^0 \oplus H^1 \oplus ...\oplus H^s \;$ we obtain
for $ \; 0 \leq s < \frac{n}{2} \;$  that

\[
L^{2} (V_s) = H^{0}(q) \oplus H^{1}(q) \oplus \cdots \oplus H^{s}(q) \hspace{2cm}
\mbox{(orthogonal sum)}
\]

\noindent
where $\; H^0(q) \;$ is the one dimensional trivial representation and $\; H^{t}(q) = (R(q)_{s-1} \circ \cdots \circ R(q)_{t}) \; (Ker
R(q)_{t-1}^{*}), \; (1 \leq t \leq s), \; $ with  dimension $\; \left(
\begin{array}{c} n \\ t
\end{array}\right)_{q}-\left( \begin{array}{c} n \\ t-1
\end{array}\right)_{q}.\;$

\vspace{12pt}

As in theorem 1, these subspaces afford all the irreducible components of the multiplicity-free natural 
representation of $\; GL_n(q) \;$   in $\; L^2(V_s),\; $  
because they have different dimensions and their number equals the
 number of orbits of   $\; GL_n(q) \;$   in     $V_s \times V_s. \;$  

\bigskip

Next we proceed to construct  the spherical function $ \;
\Phi(q)_{t} \;$ associated to  $\;  H^{t}(q).\;$ We have to find a set of
functions $\; h(q)_t \in Ker R(q)_{t-1}^{*} \;$ such that

\[
\Phi(q)_{t}(W) = \sum_{Y < W } h(q)_t(Y) =\sum_{k=k_{0}}^{min(j,t)} \left( \sum_{Y < W, \ell_q(Y,W_{0}) = k} h(q)_t(Y) \right)
\]

\noindent
where  $\;d_q(W_{0}, W) = j \; $ and $ \;  \ell_q (Y, W_{0}) = k_{0}
\; $ where $ \; k_{0} = 0 \; $ if $ \; j \leq s-t \; $ or $ \; k_{0}
= t+j-s \; $ if $ \; j \geq s-t. \; $   We require that the function
$ \; h(q)_t $ satisfies the same kind of conditions we imposed on $
\; h_t, \; $ obtaining the following expression:

\[
\Phi(q)_t(W) =  \sum_{k=k_{0}}^{\min(j,t)} |A(q)_{t}^{k,j}| \alpha(q)^{k}_{t}
\]

where $ \; \alpha(q)^{k}_{t} = h(q)_t(Y)  \; $ if $ \; \ell_q(Y, W_{0}) = k,$
and

\[
\alpha(q)_{t}^{0} = \left( \begin{array}{c} s \\ t \end{array} \right)^{-1}_q
\; \; \; \mbox{and} \; \; \;
\alpha(q)_{t}^{k} = \frac{C(q)_{t}^{0,0}C(q)_{t}^{1,1} \cdots C(q)_{t}^{k-1,k-1}}
{C(q)_{t}^{0,1}C(q)_{t}^{1,2}\cdots C(q)_{t}^{k-1,k}}\alpha(q)_{t}^{0}.
\]

Using the combinatorial results of proposition 4 we get finally the
expression for the spherical function  $ \; \Phi(q)_{t} \;$
mentioned above.

\bigskip

Next we get the eigenvalues $\; \lambda(q)_{t} \;$ of the average operator $\; M(q)_1^s \;$
 by computing $\; M(q)_1^s (\Phi(q)_t)(W_0) \;$

\[
\lambda(q)_{t} =
q(n-s)_q \left((s-t)_{q} - \frac{(s-t+1)_{q} t_q }{q (n-s)_q} \right),  \; 0 \leq t \leq s
\]

We need that $ \; \lambda(q)_{t} \neq -(n-s)_q \;$. If we compute
for which values of $\;t \;$ the equality holds,  we obtain that  $
\; \lambda(q)_{t} = -(n-s)_q \;$ if and only if $\; t = n-s. \;$
Since by hiphotesis $\; t+s < n, \;$ we conclude that $ \;
\lambda(q)_{t} \neq -(n-s)_q; \; 0 \leq t \leq s \; $ and $ \;
R(q)_{s}^{*} \circ R(q)_{s} \; $ is an automorphism. Then $ \;
R(q)_{s} \; $ is injective  for all $ \; 0 \leq s < \frac {n}{2} \; $ and

\[
L^{2} (V_{s}) \cong (Im R(q)_{0}) \oplus (Im R(q)_{0})^{\perp } \oplus
\cdots  \oplus (Im R(q)_{s-1})^{\perp}\;\;\;  (orthogonal \; sums)
\]

In order to prove  b) we follow the same procedure used in a). We
get:

\[
\varphi(q)_t(W') =  \sum_{k=k_{0}}^{\min(j,t)} |B(q)_{t}^{k,j}| \beta(q)^{k}_{t}
\]

where

\[
\beta(q)_{t}^{0} = \left( \begin{array}{c} s \\ s-t \end{array} \right)^{-1}_q
\; \; \; \mbox{and} \; \; \;
\beta(q)_{t}^{k} = \frac{D(q)_{t}^{0,0}D(q)_{t}^{1,1} \cdots D(q)_{t}^{k-1,k-1}}
{D(q)_{t}^{0,1}D(q)_{t}^{1,2}\cdots D(q)_{t}^{k-1,k}}\beta(q)_{t}^{0}
\]

\noindent
and by applying proposition 4 we prove b).
\begin{flushright}     
{$\square  $}

\end{flushright}

\vspace{12pt}

The following definition corresponds to the operator $q-$analogue to
$\; {\cal C}_{s}^{*} \;$ defined in theorem 3 of the last section.

\vspace{12pt}

\begin{hef}
 Let $\; 0 \leq s \leq \frac{n}{2}.\;$
We define $ \; {\cal C}(q)_{s}^{*} : L^{2} (V_{n-s})
\rightarrow L^{2} (V_{s}), \; $ by $ \;
{\cal C}(q)_{s}^{*}(f)(W) = \sum_{W'\oplus W = V}f(W'), \; $   for   $f \in L^{2} (V_{n-s}), W \in V_s. \; $
\end{hef}
\vspace{12pt}

\begin{pro}
\begin{itemize}
\item[a)] ${\cal C}(q)_{s}^{*} \;$ is a non trivial intertwining operator between
$\; (L^{2} (V_{n-s}),\tau') \;$ and $\; (L^{2}(V_{s}),\tau), \;$

\item[b)] ${\cal C}(q)_{s}^{*}(KerR(q)_{n-s}) = KerR(q)_{s-1}^{*}, $

\item[c)]$\;{\cal C}(q)_{s}^{*} \circ R(q)_{n-s}^{*})=
 q^{((n-s)-(s-1))}(R(q)_{s-1} \circ {\cal C}(q)_{s-1}^{*}). \;$
\end{itemize}
\end{pro}

{\flushleft{\bf Proof:}}
\begin{itemize}
\item[a)]
Let $\; g \in GL_n(q), \; f \in L^{2}(V_{n-s}) \;$ and $\; W \in V_s. \; $

\bigskip

We remark  that $\; ({\cal C}(q)_{s}^{*} \circ \tau'_g)(f)(W) =  \sum_{W'\oplus W = V}f(g^{-1}W') \;$
may be expressed as

\bigskip

$$ \sum_{W'\oplus W = V}f(g^{-1}W') = \sum_{gW'\oplus W = V}f(W'). $$

\bigskip

 It is easy to prove that the sets $\; A =\{W' \in V_S : gW' \oplus W = V \} \;$ and
$\; B =\{W' \in V_s : W' \oplus g^{-1}W = V \} \;$ are equal. If $\; g  \;$ is an automorphism
we have that

$$\sum_{W'\oplus g^{-1}W = V}f(W') = \sum_{gW' \oplus W = V}f(W'), $$

from which we obtain a).

\item[b)]
Let $\; h \in (KerR(q)_{n-s}), Z \in V_{s-1}. \;$ We have to prove that for each $\; Z \in V_{s-1}, \;$

$$\sum_{Z<Y}\sum_{Y' \oplus Y=V}h(Y') = 0 $$

is fulfilled.

\bigskip

We will prove first
that the condition $\; Z<Y,\;\; Y' \oplus Y = V \;$ is equivalent to the
condition $\; Y' \in V_{n-s}, dim(Y'\cap Z) = 0. \;$ If we let  
          $\; Y' \;$   be a subspace
of $\; V_{n-s} \;$ such that  $\; Y' \oplus Y = V \;$ for some $\; Y\;$ which 
     contains
 $\;Z \;$ then $\; dim(Y'\cap Z )= 0. \;$

\bigskip

Conversely, if $\; Y' \in V_{n-s} \;$ and $\; dim(Y'\cap Z) = 0  \;$ then we can find one dimensional
subspaces $\;<w> \;$ such that $\; V = Y' \oplus Z \oplus <w>. \;$ So, for each $\; <w> \;$
we find $\; Y = Z \oplus <w> \;$ such that $\; Z<Y \;$ and $\; Y \oplus Y' = V \;$.

\bigskip

Moreover, applying Lemma 1 we compute:

$$ |Y: Z<Y ||Y': Y \oplus Y' = V | = q^{s(n-s)}(n-(s-1))_q, $$

$$ |Y'\in V_{n-s}: dim(Y'\cap Z) = 0 | = q^{(s-1)(n-s)}\left( \begin{array}{c} n-(s-1) \\ n-s \end{array} \right)_q,
 $$

and

$$ |Y: Z<Y,\;\; Y \oplus Y' = V | = q^{n-s} . $$

Therefore we get:

\[
\sum_{Z<Y}\sum_{Y' \oplus Y=V}h(y) = q^{n-s}\sum_{ dim(Y'\cap Z) = 0 }h(Y')
\]

On the other hand,in a similar way we obtain

\[
\sum_{Z'\oplus Z = V}\sum_{Y'<Z'}h(Y') = \sum_{dim(Y'\cap Z)=0}|Z':Y'<Z',Z' \oplus Z = V |h(Y'),
\]

and we compute that:
\[
|Z': Z' \oplus Z = V||Y'<Z'| = q^{(s-1)(n-(s-1)))}\left( \begin{array}{c} n-(s-1) \\ n-s \end{array} \right)_q,
\]

\[
|Y' \in V_(n-s): dim(Y'\cap Z) = 0 | = q^{(s-1)(n-s)}
\left( \begin{array}{c} n-(s-1) \\ n-s \end{array} \right)_{q},
\]

\noindent
and
\[
|Z': Y'<Z',Z' \oplus Z = V |= q^{(s-1)},
\]

obtaining

\[
\sum_{Z'\oplus Z = V}\sum_{Y'<Z'}h(Y') = q^{s-1}\sum_{dim(Y'\cap Z)= 0}h(Y').
\]

So , for each $\;Z \in V_{s-1} \;$ we have that

\[
q^{n-s}\sum_{Z'\oplus Z = V}\sum_{Y'<Z'}h(Y') = q^{s-1}\sum_{Z<Y}\sum_{Y' \oplus Y = V}h(Y').
\]

But due to the fact that $\; h \in Ker(q)_{n-s} \;$ we have that $\; \sum_{Y'<Z'}h(Y') = 0 \;$ for each $\;Z' \in V_{n-(s-1)}. \;$

 Therefore the left side of the last equation is equal to zero
and then $\; \sum_{Z<Y} \sum_{Y' \oplus Y = V}h(Y') = 0. \;$

\item[c)]
Let $\; g \in L^2(V_{n-(s-1)}) \;$ and $\; Y \in L^2(V_{s}). \;$ We have that

 \[
({\cal C}(q)_{s}^{*} \circ R(q)_{n-s}^{*})(g)(Y)= \sum_{ Y' \oplus Y = V}( \sum_{Y'<Z'}g(Z')).
\]

First,  we note that:

\[
 \sum_{Y' \oplus Y = V}( \sum_{Y'<Z'}g(Z')) =  \sum_{dim(Z'\cap Y) = 1}K_{Z'}g(Z').
\]

where $\; K_{Z'} = |Y': Y'<Z', dim(Y'\cap Y)=0,dim(Z'\cap Y) = 1|. \;$

\bigskip

  Indeed, 
 if $\;Z' \;$ is a $\; n-(s-1) \;$-subspace of $\;V\;$ and
it satisfies that $\; Y'<Z'\;$ for  some supplementary subspace$\;
Y'\;$ of $\; Y\;$ then we can find a one dimensional subspace $\; <
u > \;$ such that $\; Z' = Y' \oplus < u >. \;$ Since $\; V = Y
\oplus Y' \;$ there exist unique vectors $\; u_1 \in Y \;$ and
$\;u_2 \in Y', \;$ such that $\; u = u_1 + u_2. \;$So we have that $\;Z' \cap Y = < u_1 > \;$and
therefore $\; dim(Z' \cap Y)=1. \;$

\bigskip

Conversely, if $\; Z'\in
V_{n-(s-1)} \;$ and $\; Z' \cap Y = < u >, \;$ then we can find an
$\;(n-s)-\;$ dimensional subspace $\; Y' \;$ such that $\; Z' = < u >
\oplus Y'.\;$

So, if there is a vector $\;w \in Y' \cap Y \;$ such
that $\; w \neq \vec{0} \;$ then  the subspace $\; < u
> \oplus < w > \;$ should be a two dimensional subspace of $\; Z' \cap Y
\;$ and this would contradict the condition $\;dim(Z' \cap Y)=1. \;$

Therefore, since $\; dim (Y' \cap Y )= 0 \;$ and $\; dimY + dimY'= n
\;$ we get that $\; Y' \;$ is a supplementary subset of $\; Y \;$
for which we have that $\; Z' < Y'. \;$

\bigskip

Moreover, applying Lemma 1 we compute:

\[
\begin{array}{ll}
|Y' \in V_{n-s}: Y' \oplus Y = V|&|Z' \in V_{n-(s-1)}: Y'<Z', Y' \oplus Y = V| \\
       $ $ & $ $                          \\
$ $ &= q^{s(n-s)} \left( \begin{array}{c} s \\ 1 \end{array} \right)_{q},
\end{array}
\]

\[
|Z'\in V_{n-(s-1)}: dim(Z'\cap Y) = 1| = q^{(s-1)(n-s)}\left( \begin{array}{c} s \\ 1 \end{array} \right)_{q},
\]

and

\[
|Y' \in V_{n-s}: Y'<Z',dim(Z'\cap Y) = 1,Y' \oplus Y = V|= q^{1(n-s)}\left( \begin{array}{c}n-s \\ n-s \end{array} \right)_{q}.
\]

In this way we get

\[
({\cal C}(q)_{s}^{*} \circ R(q)_{n-s}^{*})(g)(Y)= q^{n-s}\sum_{ dim(Z'\cap Y) = 1}g(Z').
\]

On the other hand, we note that $$\; (R(q)_{s-1} \circ {\cal C}(q)_{s-1}^{*})(g)(Y) = \sum_{Z<Y }(\sum_{Z' \oplus Z = V}g(Z')). \;$$

We claim that

\[
\sum_{Z<Y}(\sum_{Z' \oplus Z = V}g(Z')) = \sum_{dim(Z'\cap Y) = 1} J_{Z} g(Z').
\]

where $\; J_{Z}= | Z \in V_{s-1} : Z \oplus Z' = V, Z<Y,dim(Z'\cap Y) = 1|. \;$

\bigskip

If we consider an $n-(s-1)- $dimensional subspace $\; Z', \;$ such that $\; Z' \oplus Z = V \;$ for some subspace   $\; Z, \;$ which is  a subspace of $\;Y, \;$ then  we can find a vector $\; v \; $ such that $\;Y = Z \oplus < v >. \;$

 Then there are unique vectors $\; z_1 \in Z\;$ and $\; z_2 \in Z' \;$such that $\; v = z_1 + z_2 \;$and  since we have found a non zero  vector$\; z_2 \in (Y \cap Z') \;$and $\; dimZ = dimY -1 \;$ then the dimension of the subspace $\; Y \cap Z'\;$ must be one. In the same way if $\; Z' \in V_{n-(s-1)}\;$ and $\; dim(Z' \cap Y)=1 \;$ there exists an $\;(s-1)- \;$dimensional subspace $\; Z \;$such that $\; Y = (y \cap Z') \oplus Z \;$ and $\; dim(Z \cap Z') = 0. \;$ Then we have found an $\;(s-1)- \;$dimensional subspace $\; Z \;$of $\; Y \;$ such that $\; Z' \oplus Z = V.\;$   

\bigskip

\noindent
Moreover we compute:

\[
\begin{array}{l l}
|Z \in V_{s-1}: Z<Y | & | Z'\in V_{n-(s-1)}: Z' \oplus Z = V, Z<Y | =\\
$  $  & $   $  \\
 $  $ & = \left( \begin{array}{c} s \\s-1 \end{array} \right)_{q}q^{(s-1)(n-(s-1))}\left( \begin{array}{c} n-(s-1)  \\ n-(s-1) \end{array} \right)_{q} =   \\
$   $ & $   $  \\
 $    $  & = q^{(s-1)(n-(s-1))}\left( \begin{array}{c} s  \\ 1 \end{array} \right)_{q},
\end{array}
\]

\bigskip

\[
| Z'\in V_{n-(s-1)}: dim(Z'\cap Y) = 1| =
 q^{(s-1)(n-s)}\left( \begin{array}{c} s \\ 1\end{array} \right)_{q},
\]

\bigskip

and

\[
|Z \in V_{s-1} : Z \oplus Z' = V, Z<Y, dim(Z' \cap Y)=1| = q^{1(s-1)}\left( \begin{array}{c} s-1  \\ s-1 \end{array} \right)_{q},
\]

\noindent
getting in this way that

\bigskip

$$ (R(q)_{s-1} \circ {\cal C}(q)_{s-1}^{*})(g)(Y) = q^{s-1}\sum_{ dim(Z'\cap Y) = 1} g(Z').$$

\bigskip

Therefore we have

$$ \frac{1}{q^{n-s}}({\cal C}(q)_{s}^{*} \circ R(q)_{n-s}^{*})= \frac{1}{q^{s-1}} (R(q)_{s-1} \circ {\cal C}(q)_{s-1}^{*})$$

from which we get  item c).
\end{itemize}
\begin{flushright}     
{$\square  $}

\end{flushright} 
\vspace{12pt}

To prove that  $\; {\cal C}(q)_{s}^{*}\;$ is an isomorphism between the unitary natural representations $\;(L^{2} (V_{n-s}),\tau') \;$ and $\;(L^{2}(V_{s}),\tau) \;$ for  $\;0 \leq s \leq \frac{n}{2}, \;$ we need   the following lema. 
 
\bigskip

\begin{lema} The number of $\; (n-s)-
\;$dimensional subspaces $\; W' \;$ such that $\; W' \oplus W_{0} =
V \;$ and $\; d( W', W_{0}')= j\;$is given by
\[
|N_j|=\left( \begin{array}{c} n-s  \\ n-s-j \end{array} \right)_{q}\frac{((\Pi_{i=0}^{j-1}(q^{s}-q^{i})\Pi_{i=0}^{j-1}(q^{j}-q^{i})}{(q^{j}-1)}.
\]
\end{lema}

{\flushleft{\bf Proof:}}

Let $\; U \;$ be an $\; (n-s-j)- \;$dimensional subspace of $\; W' \;$
and $\; \{ u_1,...,u_{n-s-j}\} \;$a basis of $\; U. \;$ We complete
this basis with vectors $\;\{z_1,..,z_{j}\}\;$     to obtain a basis of $ \; W'\;$. So, we have that $\; \{
e_1,..,e_{s},u_1,..,u_{n-s-j}, z_1,..,z_{j}\}\;$ is a basis of $\; V.
\;$

\bigskip

To construct subspaces $\; W', \;$ for each subspace $\; U, \;$
we have to complete $\; U \;$ with vectors $\; w_{k}=
\sum_{i=1}^{s}a_i^ke_i + \sum_{t=1}^{j}b_t^kz_t , \;$ where $\; v_k=
(a_1^k,..a_s^k)\neq 0 \;$ and $\; p_k= (b_1^k,..b_j^k)\neq 0 \;$ for
$\;1 \leq k \leq j.  \;$

 Since $\;dim (W' \cap W_{0}) =0   \;$
and$\;W' \cap W_{0}'= U,  \;$ we need to choose $\; \{
v_1,..,v_{j}\} \;$and $\; \{ p_1,..,p_{j}\} \;$ linear independent,
since, if these sets are linear dependent, we can find scalars $\; c_i,
d_i, 1 \leq i \leq j, \;$ not all zeros, such that $\;
\sum_{k=1}^{j}c_kv_k= 0 \;$ or $\; \sum_{k=1}^{j}d_kp_k = 0. \;$
Then $\; \sum_{k=1}^{j}a_i^kc_k= 0,\;$ for $\; 1 \leq i \leq s  \;$ or $\;
\sum_{k=1}^{j}b_t^kd_k = 0, \;$ for $\;  1 \leq t \leq s. \;$

In this way

\[
\sum_{i=0}^{s} ( \sum_{k=1}^{j}a_i^kc_k)e_i = 0
\]

or

\[
\sum_{i=0}^{s} (\sum_{k=1}^{j}b_t^kd_k)z_k = 0.
\]

Therefore

\[
\sum_{k=1}^jc_k(\sum_{i=1}^{s}a_i^ke_i ) = 0
\]

or

\[
\sum_{k=1}^jd_k(\sum_{t=1}^{j}b_t^kz_t ) = 0.
\]

So, we obtain that $\; \sum_{k=1}^j c_kw_k \in W_0'\;$ or $\;
\sum_{k=1}^j d_kw_k \in W_0\;$ and also $\;dim (W' \cap W_{0})
\neq 0   \;$ or $\; W' \cap W_0 \neq U. \;$

\bigskip

Suppose now we have chosen $\;W' = < u_1,..,u_{n-s-j}, w_1,..,w_{j}
\;$ and $\;W'' = < u_1,..,u_{n-s-j}, w_1',..,w_{j}' > \;$ and $\;
w_{k}'= \sum_{i=1}^{s}x_i^ke_i + \sum_{t=1}^{j}y_t^kz_t. \;$

We have
that $\; W' = W'' \;$ if and only if
 there are scalars $\; l_1^k,l_2^k,..,l_j^k, 1 \leq k \leq j \;$
such that

\[
( x_1^k,x_2^k,..,x_s^k) = \sum_{i=1}^jl_i^kv_i,
\]

\[
( y_1^k,y_2^k,..,y_j^k) = \sum_{i=1}^jl_i^kp_i,
\]

for $\; 1 \leq k \leq j. \;$

Also, among all choices of  $\; (v_1,v_2,..,v_j,p_1,p_2,..,p_j) \;$
there are $\; q^{j} -1 \;$ which span the same subspace. Now, since
there are $\; (q^{s}-1)..(q^{s}-q^{j-1}) \;$ ways to choose the
vectors $\; v_1,v_2,..,v_j \;$ and $\; (q^{j}-1)..(q^{j}-q^{j-1})
\;$ ways to choose the vectors  $\; p_1,p_2,..,p_j , \;$ we have

\[
\;\frac{((\Pi_{i=0}^{j-1}(q^{s}-q^{i})\Pi_{i=0}^{j-1}(q^{j}-q^{i})}{(q^{j}-1)} \;
\]

\noindent
different ways to complete the subspace $.\; U \;$ Since there are $\; \left( \begin{array}{c} n-s  \\ n-s-j \end{array} \right)_{q}\;$ ways to choose subspaces $\; U \;$
of $\; W_0' ,\;$ we finally get 
\[
|N_j|=\left( \begin{array}{c} n-s  \\ n-s-j \end{array} \right)_{q}\frac{((\Pi_{i=0}^{j-1}(q^{s}-q^{i})\Pi_{i=0}^{j-1}(q^{j}-q^{i})}{(q^{j}-1)},
\]
\begin{flushright}     
{$\square  $}
\end{flushright}

\vspace{12pt}

\bigskip

\begin{theo} Let $ \; 0 \leq s  \leq \frac{n}{2}. \; $ The unitary natural representations $\; (L^{2}(V_{n-s}),\tau') \;$ and $\;(L^{2}(V_{s}),\tau) \;$ of the finite linear group $\; GL_n(q) \;$ are isomorphic.
\end{theo}

\vspace{12pt}

{\flushleft{\bf Proof:}} Using the previous propositions and
theorems, it is enough to prove that
 ${\cal C}(q)_{s}^{*} \;$ is injective.
\bigskip

We will prove this by induction 
 on
 $\; s \;$.

\bigskip

If $\; s = 0, \;$ we have that $\; {\cal C}(q)_{0}^{*}(f)(\{ \vec{0}
\}) = f(V), \;$ for all $\; f \in  L^{2}(V_{n}), \;$ then  $\; {\cal
C}(q)_{0}^{*} \;$ is non trivial and thus injective.
\bigskip

If $\; s = 1, \;$ we have that

\bigskip

$$\; L^{2}(V_{n-1}) = R(q)_{n-1}^{*}(L^{2}(V_{n})) \perp Ker R(q)_{n-1}. \;$$

\bigskip

Let $\; f \in L^{2}(V_{n-1}), \;$ we have that $\;  f = f_1 + f_2 \;$
 where
$\; f_1 \in R(q)_{n-1}^{*}(L^{2}(V_{n})) \;$
 and
$\; f_2 \in Ker R(q)_{n-1}. \;$

 Since $\; f_1 = R(q)_{n-1}^{*}(g), \;$ for some $\; g \in L^{2}(V_{n}) \;$ we obtain that
$\; {\cal C}(q)_{1}^{*}(f_1) = {\cal C}(q)_{1}^{*}(R(q)_{n-1}^{*}(g)). \;$

 By using  the last proposition, we get

$$\; {\cal C}(q)_{1}^{*}(f_1) = q^{n-1}(R(q)_{0}^{*} \circ {\cal C}(q)_{0}^{*})(g). \;$$

\bigskip

 If $\; {\cal C}(q)_{1}^{*}(f_1)= 0 \;$ then $\; (R(q)_{0}^{*} \circ {\cal C}(q)_{0}^{*})(g) = 0. \;$

  Since
 $\; (R(q)_{0}^{*} \;$ is  injective, we deduce $\; {\cal C}(q)_{0}^{*}(g)= 0. \;$ Then $\; g = 0 \;$ and therefore $\; f_1 = 0. \;$

 In this way we have that $\;  {\cal C}(q)_{1}^{*} \;$ restricted to the subspace $\; R(q)_{n-1}^{*}(L^{2}(V_{n})) \;$ of $\; L^{2}(V_{n-1}) \;$ is injective.

\bigskip

 From the previous proposition
$\;{\cal C}(q)_{1}^{*}(f_2) \in KerR(q)_{0}^{*}.\;$  If we compute $\; {\cal C}(q)_{1}^{*}(\varphi(q)_{1})(W_{0}) \;$ where  $\; \varphi(q)_{1}\;$ is the spherical function of $\;  KerR(q)_{n-1} \;$ we obtain:

\[
{\cal C}(q)_{1}^{*}(\varphi(q)_{1})(W_{0}) = \frac{1}{q}
\]

Then $\;{\cal C}(q)_{1}^{*} \neq 0\;$ and then $\;{\cal
C}(q)_{1}^{*} \;$ is injective.

\bigskip

We suppose now that $\;{\cal C}(q)_{t}^{*} \;$ is injective for all
$\; t < s. \;$

As in the previous case we have the decomposition

\[
 L^{2}(V_{n-s}) = R(q)_{n-s}^{*}(L^{2}(V_{n-(s-1)})) \perp Ker R(q)_{n-s}.
\]
%% new  text follows    %%

We notice first that to prove that  $\;{\cal C}(q)_{s}^{*} \;$  is injective it is enough to prove that the restriction of   $\;{\cal C}(q)_{s}^{*} \;$   to the components   $ R(q)_{n-s}^{*}(L^{2}(V_{n-(s-1)}))$  and    $ Ker R(q)_{n-s}$ is injective.

Indeed,   if
$f_1 \in R(q)_{n-s}^{*}(L^{2}(V_{n-(s-1)})),  $ then   $f_1 = R(q)_{n-s}^{*}(g) $ for some   
$g \in L^{2}(V_{n-(s-1)}) $.  
 By using item c) of last proposition  we find that

\[
{\cal C}(q)_{s}^{*}(R(q)_{n-s}^{*}(g)) = \frac{q^{n-s}}{q^{s-1}}( R(q)_{s-1} \circ {\cal C}(q)_{s-1}^{*})(g).
\]
 Therefore   $ {\cal C}(q)^\ast_{s}(f_1)$ lies in $ R(q)_{s-1} ((L^{2}(V_{n-(s-1)}).$
 
 Moreover, since   $f_2 \in  Ker R(q)_{n-s} $, we have that    $ {\cal C}(q)^\ast_{s}(f_2)$  lies in    $ R(q)^\ast_{s-1}  $  since we have proved that    $ {\cal C}(q)^\ast_{s}(R(q)_{n-s}) = Ker R(q)_{s-1}^{*} $
 
 Therefore, if  we have    $ {\cal C}(q)^\ast_{s}(f) = 0, $ then writing   $\; f = f_1+f_2, \;$ where 
 $\; f_1 \in R(q)_{n-s}^{*}(L^{2}(V_{n-(s-1)})) \;$ and $\; f_2 \in Ker R(q)_{n-s},\;$  we see that necessarily   
  $ {\cal C}(q)^\ast_{s}(f_1) = 0 $ and      $ {\cal C}(q)^\ast_{s}(f_2) = 0. $  So, if if we know that the restrictions of 
  $ {\cal C}(q)^\ast_{s} $   to   $R(q)_{n-s}^{*}(L^{2}(V_{n-(s-1)})) $  and   $ Ker R(q)_{n-s}. $   are injective, we conclude that     
 $ {\cal C}(q)^\ast_{s} $  is injective.

Now, we prove that  the restriction  of   $\; {\cal C}(q)_{s}^{*} \;$
  to the subspace $\; R(q)_{n-s}^{*}(L^{2}(V_{n-(s-1)}))
\;$ is injective.

Write  $\; f_1 \in R(q)_{n-s}^{*}(L^{2}(V_{n-(s-1)})) \;$  as  $\; f_1 \in R(q)_{n-s}^{*}(L^{2}(V_{n-(s-1)})) \;$  for some
$\; g \in L^{2}(V_{n-(s-1)}). \;$
 So, if $\;{\cal C}(q)_{s}^{*}(f_1) = 0 \;$ then $\; {\cal C}(q)_{s-1}^{*}(g) = 0. \;$

 Since $\;  R(q)_{s-1} \;$ is injective we
deduce that $\; {\cal C}(q)_{s-1}^{*}(g) = 0 \;$. Therefore by
the induction hypothesis we obtain that $\; g = 0,\;$ and therefore $\;
f_1 = 0.  \;$

To conclude the proof,  we prove now that  $\;{\cal C}(q)_{s}^{*} \;$ restricted to
the irreducible subrepresentation $\:Ker R(q)_{n-s} \;$ of $\;
L^{2}(V_{s}) \;$ is injective.

 For this we consider the spherical
function $\; \varphi(q)_{s} \;$ of the irreducible representation
$\; KerR_{n-s}^{*} \;$ and
we compute $\;{\cal C}(q)_{s}^{*}(\varphi(q)_{s})(W_{0}). \;$

We obtain that

\[
\begin{array}{l l}
{\cal C}(q)_{s}^{*}(\varphi(q)_{s})(W_{0}) & = \sum_{W' \oplus W_{0}= V}(\varphi(q)_{s})(W')\\
$  $ & $  $  \\
$   $ & = \sum_{j=0}^{s}|N_{j}|(\varphi(q)_{s})(W') \\
$  $ & $   $ \\
$   $ & = 1 + \sum_{j=1}^{s}|N_{j}|(\varphi(q)_{s})(W'),
\end{array}
\]

where $\; N_{j} = \{W' : W' \oplus W_{0}= V, d(W',W_{0}')= j \}$

We recall that :
\[
\varphi(q)_{s})(W')= \sum_{k=k_{0}}^{min(j,s)}(-1)^{k}q^{\frac{k}{2}(k-2j-1)}\left( \begin{array}{c} j  \\ k \end{array} \right)_{q}\left( \begin{array}{c} s-j  \\ s-k \end{array} \right)_{q},k!_{q}\frac{(n-s-k)!_{q}}{(n-s)!_{q}},
\]

where $\; j = d(W',W_{0}') \;$ and $\; k_{0}=j.  \;$

Then for $\; 1 \leq j \leq s \;$ we have

\[
\varphi(q)_{s})(W')= (-1)^{\frac{-j}{2}(j+1)}\frac{j!_{q}(n-s-j)!_{q}}{(n-s)!_{q}}.
\]

In this way

\[
{\cal C}(q)_{s}^{*}(\varphi(q)_{s})(W_{0})=  1 + \sum_{j=1}^{s}|N_{j}|(-1)^{j}(q)^{\frac{-j}{2}(j+1)}\frac{j!_{q}(n-s-j)!_{q}}{(n-s)!_{q}}
\]
  
From lema 2 we have that

\[
|N_j|=\left( \begin{array}{c} n-s  \\ n-s-j \end{array} \right)_{q}\frac{((\Pi_{i=0}^{j-1}(q^{s}-q^{i})\Pi_{i=0}^{j-1}(q^{j}-q^{i})}{(q^{j}-1)},
\]

Replacing $\;|N_j| \;$ in the last equation we get,

\[
{\cal C}(q)_{s}^{*}(\varphi(q)_{s})(W_{0})= 1 + \sum_{j=1}^{s}(-1)^{j}q^{\frac{j(j-3)}{2}}((\Pi_{i=0}^{j-1}(q^{s}-q^{i})\Pi_{i=1}^{j-1}(q^{i}-1).
\]

If we put $\;{\cal C}(q)_{s}^{*}(\varphi(q)_{s})(W_{0})= S + T \;$
where

\[
S = 1+\sum_{j=1}^{2}(-1)^{j}q^{\frac{j(j-3)}{2}}((\Pi_{i=0}^{j-1}(q^{s}-q^{i})\Pi_{i=1}^{j-1}(q^{i}-1)
\]

and

\[
T = \sum_{j=3}^{s}(-1)^{j}q^{\frac{j(j-3)}{2}}((\Pi_{i=0}^{j-1}(q^{s}-q^{i})\Pi_{i=1}^{j-1}(q^{i}-1),
\]

we have
\[
S = 1+(q^s-1)(q^{s-1}-q^{s-2}-1).
\]

\noindent
and

\[
T = (q^s-1)\sum_{j=3}^{s}(-1)^{j}q^{\frac{j(j-3)}{2}}((\Pi_{i=0}^{j-1}(q^{s-1}-q^{i})\Pi_{i=1}^{j-1}(q^{i}-1).
\]

In this way , if  $\; {\cal C}(q)_{s}^{*}(\varphi(q)_{s})(W_{0})= 0
\;$ then $\; T = -S \;$ i.e.

\[
(q^s-1)\sum_{j=3}^{s}(-1)^{j}q^{\frac{j(j-3)}{2}}((\Pi_{i=0}^{j-1}(q^{s-1}-q^{i})\Pi_{i=1}^{j-1}(q^{i}-1)= -1-(q^s-1)(q^{s-1}-q^{s-2}-1).
\]

\noindent
 Therefore $\; (q^s-1 )  \;$ would be a factor of $\; -1-(q^s-1)(q^{s-1}-q^{s-2}-1). \;$
But

\[
\frac{-1-(q^s-1)(q^{s-1}-q^{s-2}-1)}{(q^{s}-1 ) }= -q^{s-1}+q^{s-2}+1 + \frac{-1}{(q^{s}-1) }.
\]
 
Then $\; (q^s-1 )  \;$ is not a factor of $\;
-1-(q^s-1)(q^{s-1}-q^{s-2}-1) \;$
Therefore $\; {\cal
C}(q)_{s}^{*}(\varphi(q)_{s})(W_{0})\neq 0.  \;$ . Then  $\; {\cal C}(q)_{s}^{*}
\;$ restricted to $\;  Ker R(q)_{n-s} \;$ is injective.    
  The injectivity of   $\; {\cal C}(q)_{s}^{*}
\;$ follows.

 \begin{flushright}     
{$\square  $}

\end{flushright}
\bigskip

\begin{rem}
For $\; 0 \leq s \leq \frac{n}{2},\;$  we define
 
$$  {\cal C}(q)_{n-s}^{*} : L^{2} (V_{s})
\rightarrow L^{2} (V_{n-s}), $$

by

$$ {\cal C}(q)_{n-s}^{*}(f)(W') = \sum_{W\oplus W' = V}f(W), $$
  
 for   $\; f \in L^{2} (V_{s}), W' \in V_{n-s}.\;$

By proceeding in a similar way as we did in the proof of  proposition 5 and theorem 5, we get the following
results:

\begin{itemize}
\item[a)] ${\cal C}(q)_{n-s}^{*} \;$ is a non trivial intertwining operator between
$\; (L^{2} (V_{s}),\tau) \;$ and $\; (L^{2}(V_{n-s}),\tau'), $

\item[b)] ${\cal C}(q)_{n-s}^{*}(KerR(q)_{s-1}^{*}) = KerR(q)_{n-s}, $

\item[c)]$\;  ( {\cal C}(q)_{n-s}^{*} \circ R(q)_{s-1})= q^{((n-s)-(s-1))}(R(q)_{n-s}^{*} \circ {\cal C}(q)_{n-(s-1)}^{*}). \;$

\item[d)]$\;{\cal C}(q)_{n-s}^{*} \;$ is injective. (The $q$-analogue to the set $\;N_j \;$ is the set 
 $\; S_{j} = \{W : W \oplus W_{0}'= V, d(W,W_{0})= j \}$ and $\;|N_j | = |S_j |.\;$)

\item[e)]
$$({\cal C}(q)_{n-s}^{*}\circ{\cal C}(q)_{s}^{*}) \circ R(q)_{n-s}^{*} =
R(q)_{n-s}^{*}\circ({\cal C}(q)_{n-(s-1)}^{*}\circ{\cal C}(q)_{s-1}^{*})$$

and

$$ R(q)_{s-1} \circ({\cal C}(q)_{s-1}^{*}\circ{\cal C}(q)_{n-(s-1)}^{*})  =
({\cal C}(q)_{s}^{*}\circ{\cal C}(q)_{n-s}^{*}) \circ R(q)_{s-1}. $$

\end{itemize}
\end{rem}

\vspace{2cm}

\noindent
{\em Acknowledgements}:   The author would like to thank Profs. Christian Kassel and Jorge Soto-Andrade for stimulating and useful discussions concerning her work and Prof. Rafael Benguria for his careful reading and assessment of the manuscript.

\vspace{2cm}
\begin{flushleft}

Maria Francisca Ya\~nez

mfyanez@uchile.cl

Escuela de Pregrado

Facultad de Ciencias Qu\'imicas y Farmac\'euticas

Universidad de Chile

Olivos   1007

Santiago, Chile

\end{flushleft}

\end{document}